\documentclass[10pt,twoside,a4paper]{amsart}
\usepackage{amssymb,amsmath,latexsym,xypic,amscd,amsthm}

\newtheorem{theorem}{Theorem}[section]
\newtheorem{lemma}[theorem]{Lemma}
\newtheorem{corollary}[theorem]{Corollary}

\newtheorem{lem}[theorem]{Lemma}

\newtheorem{prop}[theorem]{Proposition}

\theoremstyle{definition}
\newtheorem{definition}[theorem]{Definition}

\theoremstyle{remark}

\usepackage{graphicx}
\usepackage{psfrag}
\usepackage{epic}
\usepackage{eepic}
\input{epsf}

\newcommand{\card}[1]{| #1 |}
\newcommand{\pro}[2]{\langle #1, #2 \rangle}

\def\Z{{\mathbb Z}}
\def\Q{{\mathbb Q}}
\def\R{{\mathbb R}}

\def\V{{{\mathcal{V}}}}
\def\P{{\mathbb P}}

\def\Vs{{\tilde{V}}}

\def\MR{{M_{\R}}}
\def\NR{{N_{\R}}}

\def\lin{{\rm lin}}
\def\conv{{\rm conv}}

\def\Hom{{\rm Hom}}

\def\Vt{{\tilde{V}}}

\def\rand{\partial}
\def\randp{{\rand P}}

\numberwithin{equation}{section}

\newcommand{\Ve}{\mathcal{V}}

\newcommand{\hyper}[2]{H(#1,#2)}
\newcommand{\hyperp}[2]{H_P(#1,#2)}
\newcommand{\nv}[2]{n(#1,#2)}
\newcommand{\nf}[2]{N(#1,#2)}

\newcommand{\pair}[2]{\langle #1,#2 \rangle}

\begin{document}

\title[$\Q$-factorial Gorenstein toric Fano varieties]
{$\Q$-factorial Gorenstein toric Fano varieties\\with large Picard number}

\author[Benjamin Nill]{Benjamin Nill}
\address{Research Group Lattice Polytopes, Freie Universit\"at Berlin, Arnimallee 3, 14195 Berlin, Germany}
\email{nill@math.fu-berlin.de}

\author[Mikkel \O bro]{Mikkel \O bro}
\address{Department of Mathematical Sciences, University of Aarhus, 8000 Aarhus C, Denmark}
\email{oebro@imf.au.dk}

\begin{abstract}
In dimension $d$, $\Q$-factorial Gorenstein toric Fano varieties with Picard number $\rho_X$ 
correspond to simplicial reflexive polytopes with $\rho_X + d$ vertices. 
Casagrande showed that any $d$-dimensional simplicial reflexive polytope 
has at most $3 d$ vertices, if $d$ is even, respectively, $3d-1$, if $d$ is odd. 
Moreover, for $d$ even there is up to unimodular equivalence only one such polytope with $3 d$ vertices, corresponding 
to $(S_3)^{d/2}$ with Picard number $2d$, where $S_3$ is the blow-up of $\P^2$ at three non collinear points. In 
this paper we completely classify all $d$-dimensional simplicial reflexive polytopes having $3d-1$ vertices, 
corresponding to $d$-dimensional $\Q$-factorial Gorenstein toric Fano varieties with Picard number $2d-1$. For $d$ even, there exist three such 
varieties, with two being singular, while for $d > 1$ odd there exist precisely two, both being nonsingular toric fiber bundles over $\P^1$. 
This generalizes recent work of the second author.
\end{abstract}

\maketitle

\vspace{-1cm}

\section{Introduction}

The main motivation of this paper is to finish the classification of simplicial reflexive polytopes with the maximal number of vertices, 
pursued in \cite{Nil05,Cas06,Oeb08}. Before stating 
the main convex-geometric result, Theorem \ref{main}, we recall necessary notions. The algebro-geometric version of Theorem \ref{main} 
is given in Corollary \ref{alggeo}.

\subsection{Lattice polytopes} 
A {\em polytope} is the convex hull of finitely many points in a vector space. Given a lattice $N \cong \Z^d$, a 
polytope $P \subseteq \NR := N \otimes_\Z \R \cong \R^d$ is called {\em lattice polytope}, if all vertices of $P$ are lattice points. We denote the 
set of vertices of $P$ by $\V(P)$. In other words, 
a lattice polytope is the convex hull of finitely many lattice points. We say two lattice polytopes are {\em isomorphic} or {\em unimodularly equivalent}, 
if there is a lattice automorphism mapping one vertex set onto the other. In what follows we always assume that 
$P$ is a lattice polytope of full dimension $d$ that 
contains the origin in its interior. In this case we can define {\em the dual polytope} $P^*$. For this let us denote by $M := \Hom_\Z(N,\Z)$ the dual lattice of $N$ and by $\MR := M \otimes_\Z \R$ the dual vector space of $\NR$. Then
$$P^* := \{x \in \MR \,:\, \pro{x}{y} \leq 1 \;\forall\, y \in P\},$$
is also a $d$-dimensional polytope containing the origin in its interior, however in general it is not a lattice polytope. 

\subsection{Reflexive polytopes}
A $d$-dimensional lattice polytope $P \subseteq \NR$ with the origin in its interior 
is called {\em reflexive polytope}, if $P^*$ is also a lattice polytope. This definition was given by Batyrev \cite{Bat94} in the context 
of mirror symmetry. It is known that 
there is only a finite number of isomorphism classes of reflexive polytopes in fixed dimension $d$, and complete 
classification results exist for $d \leq 4$, see 
\cite{KS98,KS00}. The polytope $P$ is called {\em simplicial}, if each facet (i.e., $(d-1)$-dimensional face) is a simplex. 
The most interesting case of a simplicial reflexive polytope is given by a lattice polytope containing the origin in its interior, 
where the vertices of each 
facet form a lattice basis. We call such a polytope a {\em smooth Fano polytope}. 
These special reflexive polytopes were studied quite intensely, and 
by now we have complete lists for $d \leq 8$, see \cite{Oeb07}. 

\subsection{Low dimensions}
Let us look at simplicial reflexive polytopes with many vertices in low dimensions $d$. 
For $d=1$ there is only one reflexive polytope, namely $[-1,1] \subseteq \R$ (with respect to the lattice $\Z$). 
For $d=2$ there are $16$ isomorphism classes of reflexive polytopes (all necessarily simplicial). 
Only three of these (called $\Vt_2$,$E_1$,$E_2$) have $5$ vertices, and precisely one (called $V_2$) has $6$ vertices. 
$\Vt_2,\Vt_2$ are smooth Fano polytopes, while $E_1,E_2$ are not.\\

{\centerline{ 
\psfrag{a}{$\Vt_2$}
\psfrag{b}{$E_1$}
\psfrag{c}{$E_2$}
\psfrag{d}{$V_2$}
\includegraphics{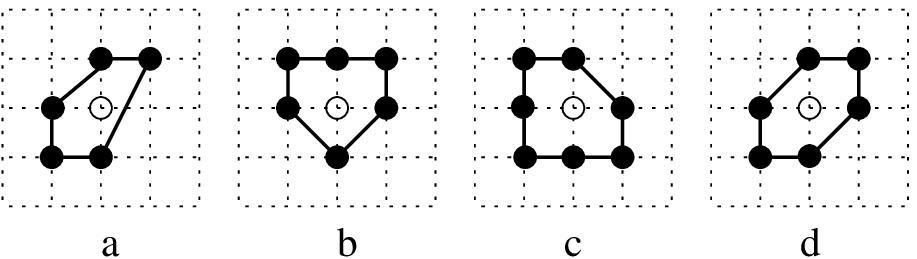}}
}

\bigskip

For $d=3$ there are $4319$ isomorphism classes of reflexive polytopes, of these are $194$ simplicial. 
There are up to isomorphisms only two three-dimensional simplicial reflexive polytopes 
having the maximal number of $8$ vertices. Both are smooth Fano polytopes 
that are bipyramids over a hexagon, we denote them by $Q_3$ and $Q'_3$:\\

\centerline{
{
\psfrag{a}{$Q_3$}
\psfrag{b}{$v$}
\psfrag{c}{$-v$}
\includegraphics[height=30mm]{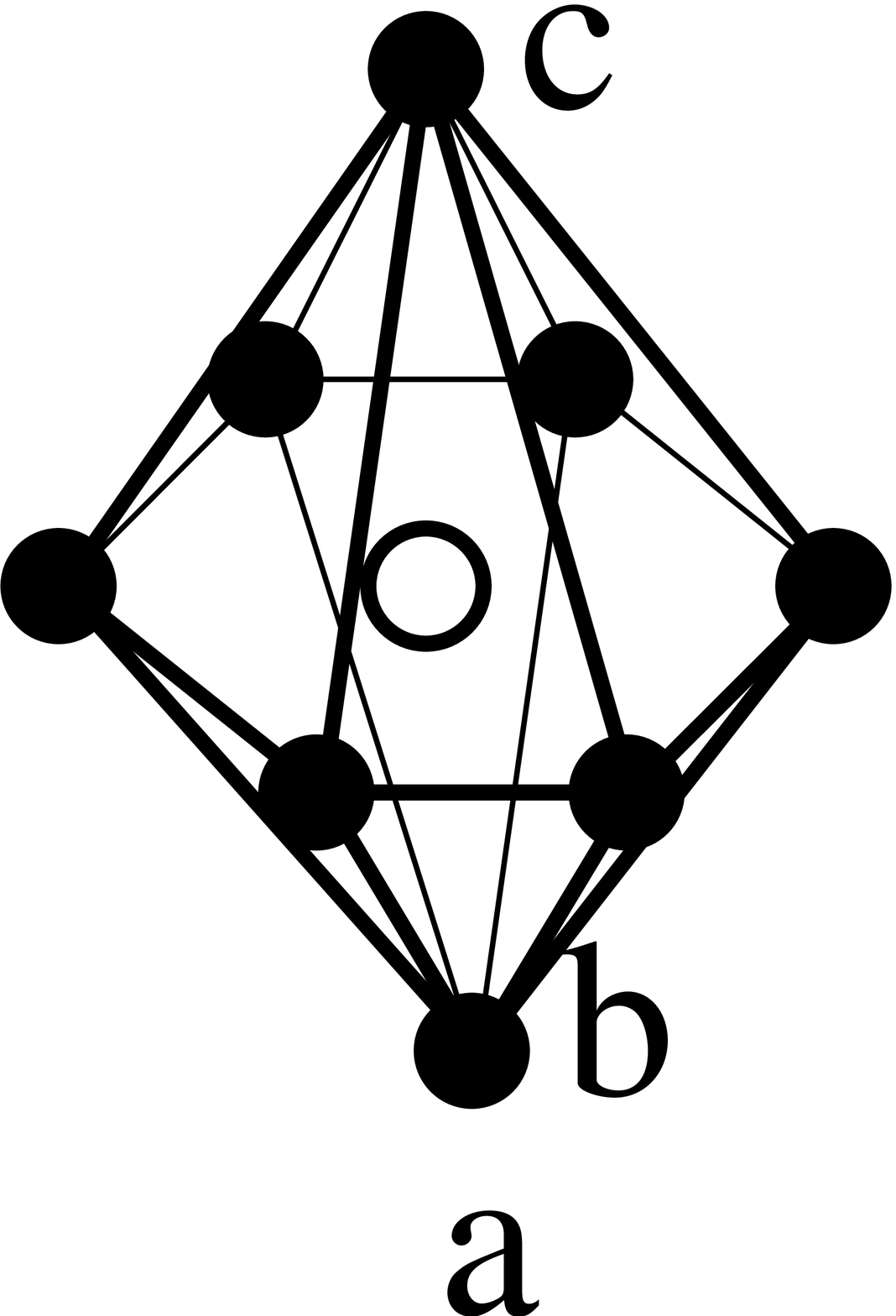}
}
{
\psfrag{a}{$Q'_3$}
\psfrag{b}{$v$}
\psfrag{c}{$v'$}
\psfrag{d}{$w$}
\includegraphics[height=30mm]{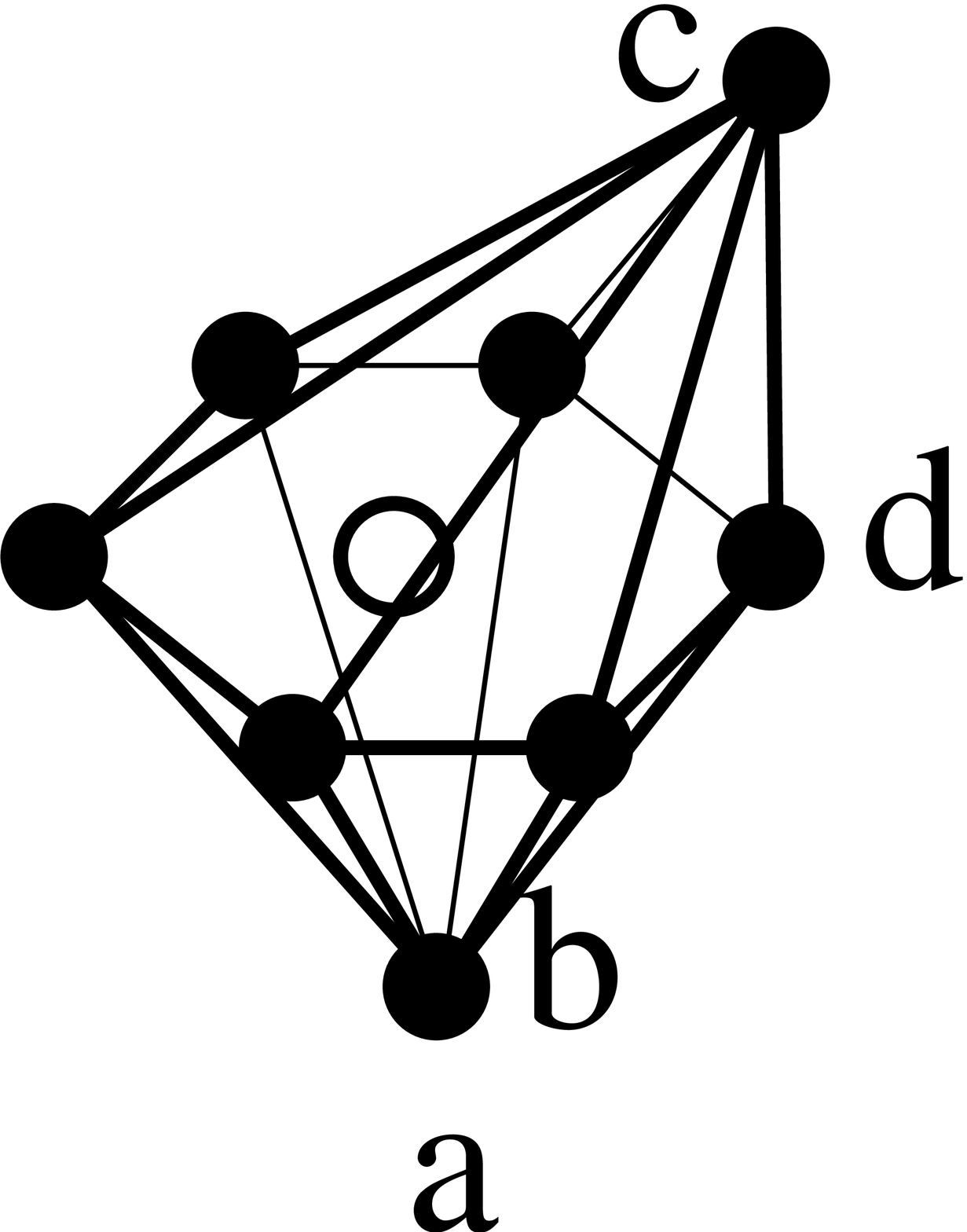}
}
}

\bigskip

While $Q_3$ is centrally symmetric, the two apexes $v,v'$ of $Q'_3$ add up to a vertex $w$ of the hexagon, i.e., $v+v'=w$.

\subsection{The main theorem}
To describe the general case, let us say a reflexive polytope $P \subseteq \NR$ {\em splits} into $P_1$ and $P_2$, 
if $P$ is the convex hull of lattice polytopes $P_1$ and $P_2$, and $N = N_1 \oplus_\Z N_2$, $P_1 \subseteq (N_1)_\R$, $P_2 \subseteq (N_2)_\R$. 
In this case, $P_1$ (respectively, $P_2$) is a reflexive polytope with respect to $N_1$ (respectively, $N_2$). 
For instance, $Q_3$ splits into $[-1,1]$ and $V_2$. 

The following long-standing conjecture on the maximal number of vertices 
was finally proven by Casagrande \cite{Cas06} in 2004 (here $|\cdot|$ denotes the cardinality): 

\begin{theorem}[Casagrande 04]
Let $P \subseteq \NR$ be a simplicial reflexive polytope of dimension $d$. 
Then 
\[\card{\V(P)} \leq \left\{\begin{array}{lcl}3d &,&d \text{ even}\\3d-1 &,&d \text{ odd}\end{array}\right.\]
If $d$ is even and $\card{\V(P)} = 3d$, then $P$ splits into $d/2$-copies of $V_2$.
\label{casa}
\end{theorem}

Note that there are by now very short proofs of this upper bound, cf. \cite{KN07,Oeb08}, see also Subsection 2.3.\\

Here is our main result, the classification of simplicial reflexive polytopes of dimension $d$ with $3d-1$ vertices.

\begin{theorem}
Let $P \subseteq \NR$ be a simplicial reflexive polytope of dimension $d \geq 3$ with $3d-1$ vertices.

If $d$ is even, then $P$ splits into $\Vt_2$ (or $E_1$, or $E_2$) and $(d-2)/2$ copies of $V_2$.

If $d$ is odd, then $P$ splits into $Q_3$ (or $Q'_3$) and $(d-3)/2$ copies of $V_2$.
\label{main}
\end{theorem}

This generalizes a recent result of the second author in \cite{Oeb08}, where this theorem was proven under the 
assumption that any lattice point on the boundary of $P$ is a vertex (for instance, if $P$ is a smooth Fano polytope). 
In this case $E_1$ and $E_2$ cannot occur, so there is only one type in Theorem \ref{main} for $d$ even.

\subsection{Algebro-geometric interpretation}

The algebro-geometric objects corresponding to reflexive polytopes $P$ are {\em Gorenstein toric Fano varieties} $X$ (i.e., 
normal complex projective varieties, where the anticanonical divisor is Cartier and ample). The relation is given 
via the toric dictionary, see \cite{Ful93}: $X$ is the toric variety associated to the fan spanned by the faces of $P$. 
For the {\em Picard number} $\rho_X$ of $X$ we have the equation $\rho_X = \card{\V(P)} - d$. For instance, $V_2$ 
corresponds to the del Pezzo surface $S_3$ with $\rho_{S_3} = 4$, which is $\P^2$ blown-up in three torus-invariant fixpoints. 
In the same way, $\Vt_2$ corresponds to the del Pezzo surface $S_2$ with $\rho_{S_2} = 3$. 
Here, $P$ is simplicial if and only if $X$ is {\em $\Q$-factorial} (i.e., for any Weil divisor some multiple is Cartier). 
Moreover, $P$ is a smooth Fano polytope if and only if $X$ is a {\em toric Fano manifold} (i.e., nonsingular). 
Since the splitting of reflexive polytopes corresponds to products of toric Fano varieties, we can reformulate Casagrande's result by saying 
that the Picard number of a $\Q$-factorial Gorenstein toric Fano variety $X$ is at most $2d$, with equality only for $d$ even and $X \cong (S_3)^{d/2}$. 
Here is the corollary of our main result:

\begin{corollary}
Let $X$ be a $\Q$-factorial Gorenstein toric Fano variety of dimension $d \geq 3$ and with Picard number $\rho_X = 2d-1$.

If $d$ is even, then $X$ is a product of $(S_3)^{(d-2)/2}$ and a (possibly singular) del Pezzo surface $S$ with $\rho_S = 3$, where there are 
three possibilities for $S$ up to isomorphisms, only one of these, namely $S_2$, is nonsingular.

If $d$ is odd, then $X$ is a product of $(S_3)^{(d-3)/2}$ and a toric Fano $3$-fold $Y$ with $\rho_Y = 5$, where 
there are two possibilities for $Y$ up to isomorphisms, namely $S_3 \times \P^1$ or a unique toric $S_3$-fiber bundle over $\P^1$.
\label{alggeo}
\end{corollary}

\subsection{Organization of this paper}

In the second section we recall preliminary results, 
namely properties of lattice points of reflexive polytopes, results about neighboring facets, and 
the notion of a special facet.

In the third section, we start the proof of the main result, which is then separated into Parts I--III, given in Sections 4--6. 
The proof is a combination of two different ideas. The first idea of the proof is the same that was successfully used in \cite{Cas06,KN07,Oeb08}: 
having a large number of vertices implies that there is a special facet from which 
nearly all vertices have integral distance two or less. Then Parts I and II can be treated using the methods developed and 
applied by the second author in \cite{Oeb08}. 
For Part III we use as a second idea the essential property of reflexive polytopes, namely their duality, 
to get restrictions on the outer normals of their facets. Then 
we can apply the strong properties of pairs of vertices of simplicial reflexive polytopes proven by the first author in \cite{Nil05}.

\subsection*{Acknowledgments}

The first author would like to thank Lutz Hille for his interest and 
discussions about the existence of exceptional sequences of line bundles on toric Fano varieties 
that initially started and motivated this project. 
Part of this research was done during a stay of the second author at the Research Group Lattice Polytopes 
at Freie Universit\"at Berlin, which is financed by Emmy Noether fellowship HA 4383/1 of the German Research Foundation (DFG).

\section{Preliminary results}

In this section we present basic results on simplicial reflexive polytopes.

\subsection{Lattice points in reflexive polytopes}

First let us recall an elementary property of reflexive polytopes, see \cite{Bat94} or \cite[Prop.1.12, Lemma 1.17]{Nil05}.

\begin{lemma}
A reflexive polytope contains no interior lattice points different from the origin. In dimension two this property is equivalent to 
the reflexivity of the polytope.
\label{reflbasic}
\end{lemma}

The following notation was introduced in \cite{Nil05}.

\begin{definition}{\rm 
Let $P$ be a polytope. We denote by $\randp$ its boundary. 

For $x,y \in \randp$, we write $x \sim y$, if $x$ and $y$ are contained in a common face (or equivalently, facet) of $P$.
}
\end{definition}

Using this relation we can describe a partial addition of lattice points in reflexive polytopes, see \cite[Prop. 4.1]{Nil05}.

\begin{lemma}
Let $P \subseteq \NR$ be a reflexive polytope, and $v, w \in \randp \cap N$.

Then $v + w \not= 0$ and $v \not\sim w$ if and only if $v + w \in \randp \cap N$.
\label{prim}
\end{lemma}

Finally, in the simplicial case there is a strong restriction on pairs of vertices, see \cite[Lemma 5.11]{Nil05}. 

\begin{lemma}
Let $P \subseteq \NR$ be a simplicial reflexive polytope. 
Let $v,w,w' \in \V(P)$ be pairwise different such that $w \not= -v \not= w'$, $v \not\sim w$ and $v \not\sim w'$. 

Then $P(v,w,w') := P \cap \lin(v,w,w')$ is a two-dimensional reflexive polytope with at least five vertices. 
\label{twolemma}
\end{lemma}

\subsection{Neighboring facets}

\smallskip
\begin{center}
{\em Throughout let $P$ be a $d \geq 2$-dimensional simplicial reflexive polytope in $\NR$.}
\end{center}
\smallskip

Let us first fix our notation. 

\begin{definition}{\rm Let $F$ be a facet of $P$.
\begin{itemize}
\item The vertices $\Ve(F)$ form a basis of $N_\R$. 
We denote by $\{u_F^v\ |\ v\in\Ve(F)\}$ the {\em dual basis} in $M_\R$, i.e., 
$\pro{u^v_F}{w} = \delta_{v,w}$ for $v,w \in \V(F)$.
\item Let $v \in \V(F)$ be a vertex of $F$. Then there is a unique facet of $P$ 
that contains all vertices of $F$ except $v$. We call this facet the {\em neighboring facet} $\nf{F}{v}$. The unique vertex 
of $\nf{F}{v}$ that is not contained in $F$ is called the {\em neighboring vertex} $\nv{F}{v}$.
\item There is a unique {\em outer normal} $u_F \in \MR$ defined by 
$\pro{u_F}{F} = 1$. The dual polytope $P^*$ has as vertices precisely the outer normals of the facets of $P$. Since $P$ is reflexive, 
the outer normal $u_F$ is a lattice point. Hence, the lattice $N$ is ``sliced'' into lattice hyperplanes
$$
\hyper{F}{i}:=\{x\in N\ |\ \pair{u_F}{x}=i\}\ ,\ i\in \Z.
$$
Let us abbreviate $\hyperp{F}{i} := \hyper{F}{i} \cap \V(P)$.
\end{itemize}
}
\end{definition}

We are going to collect restrictions on neighboring vertices and facets. 
The first result is contained in \cite[Lemma 1, Lemma 2]{Oeb08} (the point (3) follows from (1)).

\begin{lem}
\label{u_F_lemma}
Let $F$ be a facet of $P$ and $v\in \Ve(F)$. 
Let $F'$ be the neighboring facet $\nf{F}{v}$ and $v'$ the neighboring vertex $\nv{F}{v}$. 
Then
\begin{enumerate}
\item For any point $x\in N_\R$,
$$
\pair{u_{F'}}{x}=\pair{u_F}{x}+(\pair{u_{F'}}{v}-1)\pair{u_F^v}{x},
$$
where $\pair{u_{F'}}{v}-1 \leq -1$.
\label{enum1}
\item For any $x\in P$
\label{enum2}
$$
\pair{u_F}{x}-1\leq \pair{u_F^v}{x}.
$$
In case of equality, $x$ is on the facet $\nf{F}{v}$.
\item If $\nv{F}{v} \in H(F,0)$ and $\V(F)$ is a lattice basis, then $\pro{u^v_F}{\nv{F}{v}} = -1$.
\end{enumerate}
\end{lem}

For the next two results \cite[Lemma 3, Lemma 4]{Oeb08} 
compare Remark 5(2) in Section 2.3 of \cite{Deb03} and \cite[Lemma 5.5]{Nil05}.

\begin{lem}
Let $F$ be a facet, and $x$ a lattice point in $\randp \cap \hyper{F}{0}$. Then $x$ lies on a neighboring facet of $F$. In particular, 
let $v \in \hyperp{F}{0}$. Then $v$ is a neighboring vertex of $F$. Hence, there are at most $d$ vertices of $P$ in $\hyper{F}{0}$. Moreover, it holds:
\begin{enumerate}
\item For every $w\in\Ve(F)$, $v$ is equal to $\nv{F}{w}$ if and only if $\pair{u_F^w}{v}<0$. In particular, 
for every $w\in \Ve(F)$ there is at most one vertex $v \in \hyper{F}{0}$ with $\pair{u_F^w}{v}<0$.
\item If $v$ is contained in precisely one neighboring facet $\nf{F}{w}$ of $F$, then $v \not\sim w$, so $v + w \in F \cap N$. 
\end{enumerate}
\label{lemma1}
\end{lem}

\begin{lem}
\label{d-1inh0}
Let $F$ be a facet of $P$. 
Suppose there are at least $d-1$ vertices $e_1,\ldots,e_{d-1}$ in $\Ve(F)$, such that $\nv{F}{e_i}\in\hyper{F}{0}$ 
and $\pair{u_F^{e_i}}{\nv{F}{e_i}}=-1$ for every $1\leq i\leq d-1$. Then $\Ve(F)$ is a basis of the lattice $N$.
\end{lem}

The following lemma is due to the second author.

\begin{lem}
\label{basis_prop}
Assume for any facet $F$ of $P$ we have
$$
|\{\nv{F}{v}\in\hyper{F}{0}\ |\ v\in\Ve(F)\}|\geq d-1.
$$
Then there exists a facet $G$ such that $\Ve(G)$ is a $\Z$-basis of $N$.
\end{lem}

\begin{proof}
By lemma \ref{d-1inh0} we are done if there exists a facet $G$ such that the set
$$
\{v\in\Ve(G)\ |\ \nv{G}{v}\in\hyper{G}{0}\ \textrm{and}\ u_G^v(\nv{G}{v})=-1\}
$$
is of size at least $d-1$. So we suppose that no such facet exists.

Let $e_1,\ldots,e_d$ be a fixed basis of the lattice $N$ and write every vertex of $P$ in this basis. 
For every facet $F$ of $P$, we let $\det A_F$ denote the determinant of the matrix 
$$
A_F:= \left( \begin{array}{c}
v_1 \\
\vdots \\
v_d
\end{array}
\right) ,
$$
where $\Ve(F)=\{v_1,\ldots,v_d\}$. As $\det A_F$ is determined up to a sign, the number $r_F:=|\det A_F|$ is well-defined.

Now, let $F_0$ be an arbitrary facet of $P$. By our assumptions, 
there must be at least one vertex $v$ of $F_0$, such that $v'=\nv{F_0}{v}\in\hyper{F_0}{0}$ but 
$\pro{u_{F_0}^v}{v'}\neq -1$. Then $-1 < \pro{u_{F_0}^v}{v'} < 0$ by Lemma \ref{u_F_lemma}(\ref{enum2}) and 
Lemma \ref{lemma1}(1). Let $F_1$ denote the neighboring facet $\nf{F_0}{v}$. Then $r_{F_0}>r_{F_1}$.

We can proceed in this way to produce an infinite sequence of facets
$$
F_0,F_1,F_2,\ldots\ \ \ \textrm{where}\ \ \ r_{F_0}>r_{F_1}>r_{F_2}>\ldots .
$$
However, there are only finitely many facets of $P$, a contradiction.
\end{proof}

We also need \cite[Lemma 5]{Oeb08}.

\begin{lem}
\label{opp_lemma}
Let $F$ be a facet of $P$. Let $v_1,v_2\in\Ve(F)$, $v_1\neq v_2$, and set $y_1=\nv{F}{v_1}$ and $y_2=\nv{F}{v_2}$. Suppose $y_1\neq y_2$, $y_1,y_2\in\hyper{F}{0}$ and $\pair{u_F^{v_1}}{y_1}=\pair{u_F^{v_2}}{y_2}=-1$.

Then there is no vertex $x\in \Ve(P)$ in $\hyper{F}{-1}$ with $\pair{u_F^{v_1}}{x}=\pair{u_F^{v_2}}{x}=-1$.
\end{lem}

Finally, for convenience of the reader we cite Lemma 6 and Lemma 7 of \cite{Oeb08} with a 
weaker assumption. However, one checks that the proofs are precisely the same, so they are omitted.

\begin{lemma}
\label{old}
Let $F$ be a facet of $P$ such that any lattice point in $F$ is a vertex (for instance, 
$\V(F)$ is a lattice basis). If $\card{\hyperp{F}{0}} = d$, then the following holds:

\begin{enumerate}
\item $\hyperp{F}{0} = \{-y + z_y \,|\, y \in \V(F)\}$, where $z_y \in \V(F)$. Moreover, $\V(F)$ is a lattice basis.
\item If $x \in \hyperp{F}{-1}$, then $-x \in \V(F)$.
\end{enumerate}
\end{lemma}

\subsection{Special facets} 

Here we recall the crucial notion of special facets introduced by the second author in \cite{Oeb08}, which in particular 
yields a short proof 
of the upper bound in Casagrande's theorem.

The goal is to show that 
knowing the number of vertices of a $d$-dimensional simplicial reflexive polytope $P$ yields restrictions on the distribution of the vertices along the 
hyperplanes parallel to a special facet. For this, we define
\[\nu_P := \sum_{v \in \V(P)} v.\]

\begin{definition}{\rm A facet $F$ of $P$ with $\nu_P \in \R_{\geq 0} F$ is called {\em special facet}.
}
\end{definition}

Obviously, $P$ has a special facet, say $F$. Let us first deduce the following observation from the simpliciality of $P$ and Lemma \ref{lemma1}:
\begin{equation}
\label{trivial}
\card{H_P(F,1)} = d, \quad \card{H_P(F,0)} \leq d.
\end{equation}

Now, since $F$ is a special facet, we get
\begin{equation}
\label{bound}
0 \leq \pro{u_F}{\nu_P} = \sum_{v \in \V(P)} \pro{u_F}{v} = d + \sum_{i \leq -1} i \; \card{H_P(F,i)}.
\end{equation}
In particular, there are at most $d$ vertices lying in the union of hyperplanes $H(F,i)$ with $i \leq -1$. 
This yields together with Equation (\ref{trivial})
\begin{equation}
\label{total}\card{\V(P)} = d + \card{H_P(F,0)} + \sum_{i \leq -1} \card{H_P(F,i)} \leq 3d,
\end{equation}
which is the sharp upper bound in Theorem \ref{casa}. 

\section{Outline of the proof of the main theorem}

\label{outline}

\textit{For the remaining sections of this paper let $P \subseteq \NR$ be a $d \geq 3$-dimensional simplicial reflexive polytope with $3 d-1$ vertices.}
\bigskip

Let $F$ be a special facet of $P$. Taking Equations (\ref{trivial}) -- (\ref{total}) into account, we see that there are precisely three cases how 
the $3d-1$ vertices of $P$ can be distributed in the hyperplanes $\hyper{F}{i}$:\\

\begin{center}
\begin{tabular}{|c|ccc|}
\hline
 & Case A & Case B & Case C\\
\hline
$|\hyperp{F}{1}|$ & $d$ & $d$ & $d$\\
$|\hyperp{F}{0}|$ & $d$ & $d-1$ & $d$\\
$|\hyperp{F}{-1}|$ & $d-2$ & $d$ & $d-1$\\
$|\hyperp{F}{-2}|$ & $1$ & $0$ & $0$\\
\hline
\end{tabular}
\end{center}

\bigskip

Now, let us look at the lattice point $\nu_P$, which is the sum of the vertices of $P$, in the three cases A,B,C:

\begin{center}
\begin{tabular}{|c|ccc|}
\hline
 & Case A & Case B & Case C\\
\hline
$\pro{u_F}{\nu_P}$ & $0$ & $0$ & $1$\\
\hline
\end{tabular}
\end{center}

\bigskip

Hence, the definition of a special facet implies: {\em In the cases A and B the sum of all the vertices of $P$ equals the origin, 
while in case C the sum is a lattice point on the facet $F$.}\\

Now, the proof falls into Parts I--III (Sections 4--6) 
depending on whether $\nu_P = 0$, $\nu_P$ is a vertex, or otherwise. Then the main result, Theorem \ref{main}, follows directly from combining 
Propositions \ref{part1}, \ref{part2}, \ref{part3}.

\section{Part I: $\nu_P = 0$}

Here, we prove the following result:

\begin{prop}
Let $\nu_P = 0$. Then either $d$ is even and 
$P$ splits into $\frac{d-2}{2}$ copies of $V_2$ and a single copy of the polytope $E_2$, or $d$ is odd and $P$ splits into $\frac{d-3}{2}$ copies of $V_2$ 
and a single copy of the polytope $Q_3$.
\label{part1}
\end{prop}

\begin{proof}

Since $\nu_P=0$, every facet of $P$ is special. 
Thus, for any facet $F$ of $P$ we are in cases A or B, described above. In particular, 
there are at least $d-1$ vertices in $\hyper{F}{0}$, hence 
by Lemma \ref{lemma1} the assumptions of Lemma \ref{basis_prop} are satisfied, so we find a facet $F$ 
whose vertex set $\Ve(F)$ is a lattice basis of $N$. Let us denote the vertices of $\Ve(F)$ by $e_1,\ldots,e_d$.\\

\textbf{Claim:} We may assume we are in case A.
\smallskip

\begin{proof}[Proof of Claim]

Suppose not. Then there are $d$ vertices in $\hyper{F}{-1}$.
\smallskip

Let us first consider the case that $P$ contains a centrally symmetric pair of facets. 
Then from Theorem 0.1 in \cite{Nil07} one easily derives that 
either $d$ is even and $P$ splits into $\frac{d-2}{2}$ copies of $V_2$ and a single copy of the polytope $\Vs_2$, or $d$ is odd and 
$P$ splits into $\frac{d-3}{2}$ copies of $V_2$ and a single copy of the polytope $Q_3$. In the first case we have a contradiction to $\nu_P = 0$, 
while the second case is as desired.
\smallskip

Hence, we may assume that at least one of the vertices in $\hyper{F}{-1}$ is not equal to $-e_i$ for $i \in \{1, \ldots, d\}$, so this vertex 
has at least one positive $e_j$-coordinate for some $j$. Say, $w\in\hyperp{F}{-1}$ and $\pair{u_F^{e_1}}{w}>0$. 
Then $\pair{u_{\nf{F}{e_1}}}{w}<-1$ by Lemma \ref{u_F_lemma}(1), 
which implies that the vertices of $P$ are distributed in hyperplanes $\hyper{\nf{F}{e_1}}{\cdot}$ as in case A. In particular, 
$\pro{u_{\nf{F}{e_1}}}{w} = -2$. Now, it remains to show that $\Ve(\nf{F}{e_1})$ is a lattice basis.

If $\nv{F}{e_1}\in\hyper{F}{0}$, then $\pair{u_F^{e_1}}{\nv{F}{e_1}}=-1$ by Lemma \ref{u_F_lemma}(3), hence 
$\Ve(\nf{F}{e_1})$ is a lattice basis, as desired. So suppose $\nv{F}{e_1}\notin\hyper{F}{0}$, thus $\nv{F}{e_1} \in \hyperp{F}{-1}$. 
Since $\card{\hyperp{F}{0}} \geq d-1$, we have $\nv{F}{e_2},\ldots,\nv{F}{e_d}\in\hyper{F}{0}$ and they are all distinct. Furthermore, 
Lemma \ref{u_F_lemma}(3) yields
$$
\pair{u_F^{e_2}}{\nv{F}{e_2}}=\ldots=\pair{u_F^{e_d}}{\nv{F}{e_d}}=-1.
$$
By Lemma \ref{u_F_lemma}(2) we get $\pair{u_F^{e_i}}{w} \geq -2$; and moreover, 
if $\pair{u_F^{e_i}}{w}=-2$ for some $i>1$, then $w=\nv{F}{e_i} \in \hyper{F}{0}$, which is not possible. So $\pair{u_F^{e_i}}{w}\geq -1$ for $i = 2, \ldots, d$. 
Now, since 
$$\pair{u_F^{e_1}}{w}>0 \text{ and } \sum_{i=1}^d \pro{u_F^{e_i}}{w} = -1,$$ 
there are at least two indices $i\neq j$ in $\{2, \ldots, d\}$ such that $\pair{u_F^{e_i}}{w}=\pair{u_F^{e_j}}{w}=-1$. 
This is a contradiction to Lemma \ref{opp_lemma}.
\end{proof}

\pagebreak
So we may safely assume that $\Ve(F)$ is a lattice basis and there are $d$ vertices of $P$ in $\hyper{F}{0}$, 
$d-2$ in $\hyper{F}{-1}$ and a single one, say $v$, in $\hyper{F}{-2}$. If \mbox{$\pair{u_F^{e_i}}{v}>0$} for some $i$, 
then $\pair{u_{\nf{F}{e_i}}}{v}<-2$ by Lemma \ref{u_F_lemma}(1), which cannot happen. 
Furthermore $v$ cannot be equal to $-2e_i$ for some $i$. So (up to renumeration) $v=-e_1-e_2$, since $\pro{u_F}{v} = -2$.
The vertices of $P$ in $\hyper{F}{0}$ are by Lemma \ref{old}(1)
$$
\nv{F}{e_1}=-e_1+e_{i_1}\ ,\ \ldots \ ,\ \nv{F}{e_d}=-e_d+e_{i_d},
$$
for $\{i_1, \ldots, i_d\} \subseteq \{1, \ldots, d\}$.
By Lemma \ref{old}(2) the vertices in $\hyper{F}{-1}$ are
$$
-e_{j_1}\ ,\ \ldots \ ,\ -e_{j_{d-2}},
$$
for $\{j_1, \ldots, j_{d-2}\} \subseteq \{1, \ldots, d\}$. Suppose $i_1\neq 2$. Then there are two cases:
\begin{description}
\item[$i_2\neq 1$] It is proven precisely as in the proof of Case 2 of the main result in \cite{Oeb08} 
(starting from the line "Let $G = \nf{F}{e_1}$", with $j=i_1$ and $i=i_2$) 
that this case leads to a contradiction. 
\item[$i_2=1$] Consider the facet $G=\nf{F}{e_2}$. Since $v = -e_1-e_2 = -2 e_1 + \nv{F}{e_2}$, we see 
$\pair{u_G}{v}=-1$ and $\pair{u_G^{e_1}}{v}=-2$. Then by Lemma \ref{u_F_lemma}(\ref{enum2}), $v=\nv{G}{e_1}$. 
However, $-e_1+e_{i_1} \in \hyperp{G}{0}$ is also equal to $\nv{G}{e_1}$ by Lemma \ref{lemma1}(1), a contradiction.
\end{description}
So $i_1=2$. By symmetry, $i_2=1$. 
Now, as we see from Figure \ref{split}, $-e_1$ and $-e_2$ cannot be vertices (here, $\conv$ denotes the convex hull).

\begin{figure}[h]
{\centerline{ 
\psfrag{a}{$e_1$}
\psfrag{b}{$e_2$}
\psfrag{c}{$v$}
\psfrag{d}{$\nv{F}{e_2}$}
\psfrag{e}{$\nv{F}{e_1}$}
\includegraphics{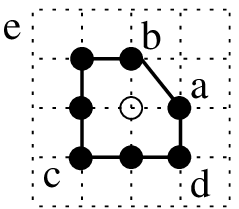}}
}
\caption{\label{split}$\conv(e_1,e_2,v,\nv{F}{e_1},\nv{F}{e_2})$}
\end{figure}

Hence,
$$
\hyperp{F}{-1}=\{-e_3,\ldots,-e_d\}.
$$
Let us consider the facet $\nf{F}{e_3}$. Since $\nv{F}{e_3} = -e_3 + e_{i_3}$, the vertices of $\nv{F}{e_3}$ form 
a lattice basis. And, because we have $v \in \hyperp{\nf{F}{e_3}}{-2}$, we are still in case A. Consequently,
$$
\hyperp{\nf{F}{e_3}}{-1}=-\Ve(\nf{F}{e_3})\setminus\{-e_1,-e_2\}.
$$
Since $-e_3+e_{i_3}\in\Ve(\nf{F}{e_3})$, the point $-e_{i_3}+e_3$ is a vertex of $P$.
\smallskip

From this, we conclude that the vertices of $P$ in $\hyper{F}{0}$ 
come in centrally symmetric pairs, so $d$ is even and $P$ splits into the claimed polytopes.
\end{proof}

\section{Part II: $\nu_P$ is a vertex of $P$}

Here, we prove the following result:

\begin{prop}
Let $\nu_P$ be a vertex of $P$. Then either $d$ is even and $P$ splits into $\frac{d-2}{2}$ copies of $V_2$ and a single copy of the polytope $\Vt_2$, 
or $d$ is odd and $P$ splits into $\frac{d-3}{2}$ copies of $V_2$ and a single copy of the polytope $Q'_3$.
\label{part2}
\end{prop}

\begin{proof}

A facet of $P$ is special if and only if it contains $\nu_P$.\\

\textbf{Claim:} Let $F$ be a special facet. Then $\Ve(F)$ is a lattice basis. 

\begin{proof}[Proof of Claim]

Since we are in case C, 
the vertices of $P$ are distributed in hyperplanes $\hyper{F}{\cdot}$ like this: $d$ in $\hyper{F}{1}$, $d$ in $\hyper{F}{0}$ and $d-1$ in $\hyper{F}{-1}$. In particular, $\nv{F}{w}\in\hyper{F}{0}$ for every $w\in\Ve(F)$.

Consider the facet $\nf{F}{w}$ for some $w\in\Ve(F)$, $w\neq \nu_P$. 
Since $\nf{F}{w}$ is also special, there are $d$ vertices in $\hyper{\nf{F}{w}}{0}$. 
So $w\in\hyper{\nf{F}{w}}{0}$ and $\nv{F}{w}\in\hyper{F}{0}$, and it follows from Lemma \ref{u_F_lemma}(1) that $\pair{u_F^w}{\nv{F}{w}}=-1$. 
This holds for all $w\in\Ve(F)$, $w\neq \nu_P$, and Lemma \ref{d-1inh0} yields that $\Ve(F)$ is a lattice basis.
\end{proof}

Let $\nu_P = e_1$. Now, the remaining proof follows precisely as in Case 1 of the proof of the 
main result in \cite{Oeb08} (starting from line "There are $d-1$ vertices in $\hyper{F}{-1}$"). 
The only difference is that in our situation one refers to points (1) or (2) in Lemma \ref{old} instead of refering to 
Lemmas 6 or 7 in \cite{Oeb08}.
\end{proof}

\section{Part III: $\nu_P \not= 0$ is not a vertex of $P$} 

Here, we prove the following result, which finishes the proof of Theorem \ref{main}:

\begin{prop}
Let $\nu_P \not= 0$, and let $\nu_P$ be not a vertex of $P$. 
Then $d$ is even and $P$ splits into $\frac{d-2}{2}$ copies of $V_2$ and a single copy of the polytope $E_1$.
\label{part3}
\end{prop}

\begin{proof}

Let $F$ be a special facet of $P$. As described in Section \ref{outline} we are in Case C, and 
$\nu_P$ is a lattice point of $F$ but not a vertex.

Let $\V(F) = \{e_1, \ldots, e_d\}$ and $\hyperp{F}{0} = \{v_1, \ldots, v_d\}$. 
By Lemma \ref{lemma1} we may assume that $v_i = \nv{F}{e_i}$ for $i = 1, \ldots, d$. In particular, Lemma \ref{lemma1}(2) 
implies:\\

\textbf{Fact~1:} For $i = 1, \ldots, d$ we have $v_i + e_i \in F$.\\

Moreover, since any neighboring vertex of $F$ is in $\hyper{F}{0}$, Lemma \ref{lemma1} implies:\\

\textbf{Fact~2:} Any lattice point in $\randp \cap \hyper{F}{0}$ is a neighboring vertex of $F$.\\

Let $G := \{x \in P \,:\, \pro{u}{x} = 1\}$. Since 
$G$ contains $d-1$ vertices, $G$ is a $(d-2)$-dimensional face of $P$. Let $\V(G) = \{b_1, \ldots, b_{d-1}\}$. 
There exist precisely two facets $G_1,G_2$ of $P$ containing $G$. We have $G = G_1 \cap G_2$. 
Let $w_1 := u_{G_1}$, $w_2 : = u_{G_2}$. The next observation is the crucial starting point of our proof.\\

\textbf{Claim 1:} By possibly interchanging $w_1$ and $w_2$ we have
\begin{enumerate}
\item[(1)] $2 w_1 + w_2 + 3 u = 0$ or
\item[(2)] $w_1+w_2+2 u = 0$.
\end{enumerate}

\bigskip

\begin{proof}[Proof of Claim 1]

By duality $w_1,w_2$ are vertices of $P^*$ joined by an edge that contains $-u$ in its relative interior. 
Let $T := \conv(w_1,w_2,u)$. By Lemma \ref{reflbasic} $T$ does not contain any lattice points different from the origin in its interior, 
thus it is a reflexive polygon. Let us look at the list of two-dimensional reflexive triangles, see, e.g., \cite[Prop.2.1]{Nil05}:\\

\centerline{
{
\includegraphics[width=100mm]{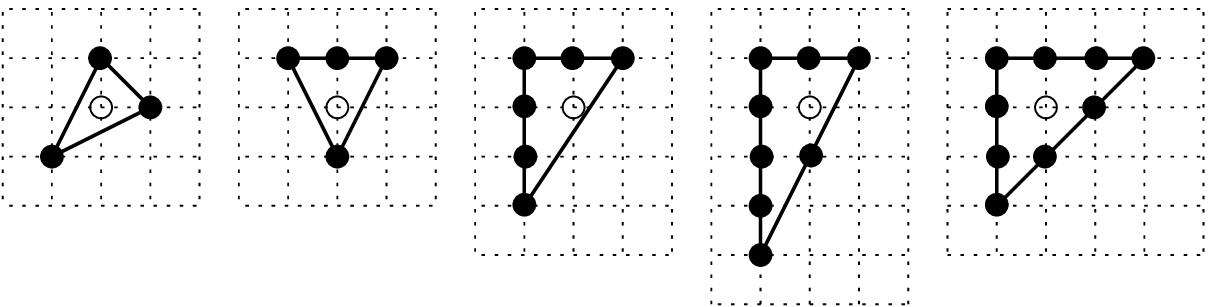}
}
}

\bigskip

From this figure we can read off all possibly occuring relations among the vertices.
\end{proof}

Now, let $x_1 \in \V(G_1)$, $x_1 \not\in G$; $x_2 \in \V(G_2)$, $x_2 \not\in G$, i.e., 
$x_1 = \nv{G_2}{x_2}$ and $x_2 = \nv{G_1}{x_1}$. 
By using $0 = \pro{2 w_1+w_2+3u}{v_i}$ in case (1) of Claim 1, respectively $0 = \pro{w_1+w_2+2u}{v_i}$ in case (2), 
we deduce:\\

\textbf{Fact~3:} Let $i \in \{1, \ldots, d\}$ such that $v_i \not\in \{x_1,x_2\}$. Then $v_i \in \hyper{G_1}{0} \cap \hyper{G_2}{0}$. 
In particular, by Lemma \ref{lemma1}, $v_i$ is a neighboring vertex of $G_1$, as well as of $G_2$.\\

\textbf{Claim 2:} $x_1$ and $x_2$ are in $\hyper{F}{0}$.

\begin{proof}[Proof of Claim 2]

Assume not. First let us suppose that $x_1 \in F$ and $x_2 \in F$. Then 
by Fact~3 and Lemma \ref{lemma1} we have $\hyperp{G_1}{0} = \{v_1, \ldots, v_d\}$. 
Moreover, Lemma \ref{lemma1} implies $x_2 = \nv{G_1}{x_1} \in \{v_1, \ldots, v_d\}$, a 
contradiction.

Since we are not going to distinguish between cases (1) and (2) for the proof of Claim~2, we may assume that 
$x_1$ is in $\hyper{F}{0}$ and $x_2$ is in $F$. Let us suppose $x_1 = v_1$. Then by Fact~3, 
$x_2, v_2, \ldots, v_d$ are the $d$ different neighboring vertices of $G_1$. By Lemma \ref{lemma1}(2) we may permute $b_2, \ldots, b_d$ 
such that $v_2 + b_2, \ldots, v_d + b_d \in G_1$. Now, by Fact~1, 
Lemma \ref{twolemma} implies for $i \in \{2, \ldots, d\}$ that $P(v_i,e_i,b_i)$ is a reflexive polygon with at least five 
vertices. Looking at the four figures in Subsection 1.3 of the introduction we find using Fact~2 that only the following two possibilities may arise:\\

\begin{figure}[h]
{\centerline{ 
\psfrag{a}{$v_i$}
\psfrag{b}{$e_i$}
\psfrag{c}{$b_i$}
\psfrag{d}{$v_i$}
\psfrag{e}{$e_i$}
\psfrag{f}{$b_i$}
\includegraphics{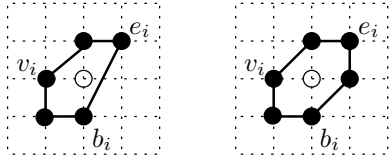}}
}
\caption{\label{picp1}Two possibilities for $P(v_i,e_i,b_i)$}
\end{figure}

\smallskip

In particular, $-b_i = v_i + e_i \in F$. Moreover, since, by Fact~3, $v_i \in \hyper{G_2}{0}$, 
we note $\pro{w_2}{e_i} = -1$. Hence, $e_i \not\in G_2$ for $i = 2, \ldots, d$. 
Therefore, $x_2 = e_1$. This implies $\pro{w_2}{e_1} = \pro{w_2}{x_2} = 1$. Since $\pro{w_2}{-b_i} = -1$, we get 
$-b_i \in \conv(e_2, \ldots, e_d)$ for $i = 2, \ldots, d$. Hence, 
$$-G = -\conv(b_2, \ldots, b_d) \subseteq \conv(e_2, \ldots, e_d).$$ 
Since, by Figure \ref{picp1}, 
$-e_2, \ldots, -e_d \in G$, this yields $G = - \conv(e_2, \ldots, e_d)$, and $\{b_2, \ldots, b_d\} = \{-e_2, \ldots, -e_d\}$. 

We conclude that for $i = 2, \ldots, d$ each vertex $v_i$ may be written as $v_i = -b_i - e_i = e_j - e_i$ for some $j \in \{2, \ldots, d\}$ 
with $j \not= i$. Hence, $\pair{u_F^{e_i}}{v_i}=-1$ for $i = 2, \ldots, d$. Now, Lemma \ref{d-1inh0} yields that $e_1, \ldots, e_d$ is a 
lattice basis. Thus, any lattice point in $F$ is a vertex, in particular, this holds for $\nu_P$, a contradiction.
\end{proof}

So Claim 2 is proven. 

Assume we are in case (1) of Claim 1. Then 
$0 = \pro{2 w_1 + w_2 + 3 u}{x_2} = 2 \pro{w_1}{x_2} + 1$, thus $\pro{w_1}{x_2} = -1/2 \not\in \Z$, a contradiction. 
Hence, we are in case (2). We may suppose $x_1 = v_1$ and $x_2 = v_2$. 
Now, Fact~3 implies
\[H_P(G_1,1) = \{x_1,b_2,\ldots, b_d\},\quad H_P(G_2,1) = \{x_2,b_2,\ldots, b_d\},\]
\[H_P(G_1,0) \supseteq \{v_3,\ldots, v_d\},\quad H_P(G_2,0) \supseteq \{v_3,\ldots, v_d\}.\]

Moreover, by $w_1+w_2+2u=0$ we get 
\[x_2 \in H_P(G_1,-1),\quad x_1 \in H_P(G_2,-1).\]

In particular, since $x_2$ is a neighboring vertex of $G_1$, however $x_1 \not\in H_P(G_1,0)$, 
Lemma \ref{lemma1} implies the following observation (the same argument holds for $G_2$):\\

\textbf{Fact~4:} $\card{H_P(G_1,0)} \leq d-1$ and $\card{H_P(G_2,0)} \leq d-1$.\\

It is our next goal to determine on which slices with respect to $G_1$ and $G_2$ the vertices of $F$ lie. 
For this we need a preliminary result.\\

\textbf{Claim 3:} $\{y \in P \,:\, \pro{w_1}{y} = -1\}$ is not a face of $P$. The same statement holds for $w_2$.

\begin{proof}[Proof of Claim 3]

Assume the claim is wrong for $w_1$. Hence, $|\hyperp{G_1}{-1}| \leq d$, and Fact~4 yields 
\[(\card{H_P(G_1,1)}, \card{H_P(G_1,0)}, \card{H_P(G_1,-1)}) = (d,d-1,d).\]
In this case $\pro{w_1}{\nu_P} = 0$, thus $\nu_P \in F \cap H(G_1,0)$. 
Now, since $v_3, \ldots, v_d \in H_P(G_1,0)$, there exists $i \in \{1, \ldots, d\}$ such that $e_i \in H_P(G_1,0)$, 
while $e_j \in H_P(G_1,-1)$ for $j \not= i$. Hence, $\nu_P \in F \cap H(G_1,0)$ implies $\nu_P = e_i \in \V(P)$, a contradiction.
\end{proof}

Since $w_1+w_2+2u=0$, we have $\pro{w_1+w_2}{e_i} = -2$ for $i = 1, \ldots, d$. Therefore, Claim 3 implies the existence of 
\[r \in \{1, \ldots, d\} \,:\, \pro{w_1}{e_r} = -2, \quad \pro{w_2}{e_r} = 0,\]
\[s \in \{1, \ldots, d\} \,:\, \pro{w_1}{e_s} = 0, \quad \pro{w_2}{e_s} = -2.\]
Moreover, since $v_2 = x_2 = \nv{G_1}{x_1}$, we get by Fact~4 and Lemma \ref{lemma1} 
\[\{e_s, v_3, \ldots, v_d\} = \{\nv{G_1}{b_2}, \ldots, \nv{G_1}{b_d}\}.\] 
We may permute $b_2, \ldots, b_d$ such that $e_s = \nv{G_1}{b_2}$ and $v_i = \nv{G_1}{b_i}$ for $i = 3, \ldots, d$; 
moreover, by Lemma \ref{lemma1}(2) we have $v_i \not\sim b_i$. Hence by Fact~1 we may apply 
for $i = 3, \ldots, d$ Lemma \ref{twolemma} to $v_i,e_i,b_i$, 
and deduce as in the proof of Claim~2 the following result:\\

\textbf{Fact~5:} For $i = 3, \ldots, d$ the polygon $P(v_i,e_i,b_i)$ looks as in Figure \ref{picp1}.\\

In particular, $\pro{w_1}{e_i} = \pro{w_2}{e_i} = - 1$ for $i = 3, \ldots, d$. 
Thus, $\{r,s\}=\{1,2\}$. Since by Fact~1, $e_1 \not\sim v_1 = x_1$, however $e_s = \nv{G_1}{b_2} \sim x_1 \in \nf{G_1}{b_2}$, we get $s \not= 1$. Hence, 
$r=1,s=2$. Let us sum up what we just proved:

\[H_P(G_1,0) = \{e_2,v_3, \ldots, v_d\}, \quad H_P(G_2,0) = \{e_1,v_3, \ldots, v_d\}.\]
\[H_P(G_1,-1) = \{x_2,e_3,\ldots, e_d\}, \quad H_P(G_2,-1) = \{x_1,e_3,\ldots, e_d\},\]
\[H_P(G_1,-2) = \{e_1\}, \quad H_P(G_2,-2) = \{e_2\}.\]

\bigskip

Now, since $e_2 = \nv{G_1}{b_2}$ we have by Lemma \ref{lemma1}(2) that $e_2 + b_2 \in G_1 \cap H(F,0)$. Thus, Fact~2 implies 
$e_2 + b_2 = x_1 = v_1 \in \V(P)$. By Fact~1 we may again apply Lemma \ref{twolemma} to $e_2,v_2,b_2$ to deduce that 
$P(v_2,e_2,b_2)$ is a reflexive polygon that has to look like one of the following two reflexive polygons (use $v_2 \sim b_2$):

\begin{figure}[h]
{\centerline{ 
\psfrag{a}{$z$}
\psfrag{b}{$e_2$}
\psfrag{c}{$v_1$}
\psfrag{d}{$b_2$}
\psfrag{e}{$v_2$}
\psfrag{f}{$e_2$}
\psfrag{g}{$v_2$}
\psfrag{h}{$b_2$}
\psfrag{i}{$z$}
\psfrag{j}{$v_1$}
\includegraphics{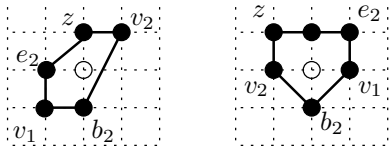}}
}
\caption{\label{picp2} Two possibilities for $P(v_2,e_2,b_2)$}
\end{figure}

\pagebreak

\textbf{Claim 4:} In Figure \ref{picp2} only the right possibility occurs, moreover, $z=e_1$. 
In particular, $P(v_2,e_2,b_2) \cong E_1$.\\

\begin{proof}[Proof of Claim 4]

Assume $P(v_2,e_2,b_2)$ is given by the left reflexive polygon in Figure \ref{picp2}. 
Since $e_1 \in H(G_2,0)$, Lemma \ref{lemma1}(2) implies that $e_1 = \nv{G_2}{b_j}$ for some $j \in \{2, \ldots, d\}$, 
even more, $e_1 \not\sim b_j$ and $e_1 \not= -b_j$. 
If $j \in \{3, \ldots, d\}$, then by Lemma~\ref{twolemma}, $e_1 \in P(b_j,v_j,e_1)=P(v_j,e_j,b_j)$, 
which implies by Figure~\ref{picp1} that $e_1 = -b_j$, a contradiction. 
Hence $j=2$. 
Therefore, Lemma \ref{twolemma} implies $e_1 \in P(b_2,e_2,e_1)=P(v_2,e_2,b_2)$. Figure \ref{picp2} yields $e_1=z=-b_2$, a contradiction.

Finally, note that in the right reflexive polygon $z \in \nf{G_1}{-2}$, and therefore $z = e_1$.
\end{proof}

By Fact~5 and Figure \ref{picp1} we know $v_i + e_i = - b_i \in F$ for $i = 3, \ldots, d$.\\

\textbf{Claim 5:} Let $i \in \{3, \ldots, d\}$. Then $-b_i \in \conv(e_3, \ldots, e_d)$.

\begin{proof}[Proof of Claim 5]

Assume not. For $j \in \{1, \ldots, d\}$ let $F_j := \nf{F}{e_j}$ and $u_j := u_{F_j}$. 
For $j = 3, \ldots, d$ we deduce from Figure \ref{picp1} that $\pm e_j \in P$ and, of course, $\pm e_j \not\in F_j$. 
This implies
\begin{equation}
\pro{u_j}{e_j} = 0 \text{ for } j = 3, \ldots, d.\label{eq2}
\end{equation}
Let $j \in \{3, \ldots, d\}$. Assume $-b_i \not\in F_j$. Then, since also $b_i \not\in F_j$, we get 
\mbox{$\pro{u_j}{-b_i} = 0$}. Now, since $-b_i \in F$, Equation (\ref{eq2}) yields $-b_i = e_j$, a contradiction to our assumption. 
Therefore, $-b_i \in F_j$ for all $j \in \{3, \ldots, d\}$. 
Hence, $-b_i \in \conv(e_1, e_2)$. Now, looking at Figure \ref{picp2} yields $b_i=b_2$, a contradiction.
\end{proof}

By Figures \ref{picp1} and \ref{picp2} we have 
\[-\conv(\frac{e_1+e_2}{2}, e_3, \ldots, e_d) \subseteq G = \conv(b_2, b_3, \ldots, b_d).\]
Now, Claim 5 shows that equality holds. Moreover, we get
\[\{-e_3, \ldots, -e_d\} = \{b_3, \ldots, b_d\}.\]
Hence, there exists a permutation $\sigma$ on $\{3, \ldots, d\}$ satisfying
\[e_{\sigma(i)} = -b_i \text{ and } \sigma(i) \not= i\quad  \text{ for } i = 3, \ldots, d.\]
By Fact~1 and Figure \ref{picp1}, $v_{\sigma(i)} \not\sim e_{\sigma(i)} = - b_i \not\sim -e_i = b_{\sigma^{-1}(i)}$, 
thus by Lemma \ref{twolemma} we have $v_{\sigma(i)} \in P(e_{\sigma(i)},v_{\sigma(i)},b_{\sigma^{-1}(i)})=P(v_i,e_i,b_i)$. 
Hence, we see that in Figure~\ref{picp1} the first possibility cannot occur, so $P(v_i,e_i,b_i) \cong V_2$, and we have 
\[v_{\sigma(i)} = -v_i \text{ for } i = 3, \ldots, d.\]
Therefore, $\sigma$ is a fix-point-free involution, thus, $\sigma$ is a product of disjoint transpositions. 
In particular, $d$ is even. It remains to show the following statement.\\

\textbf{Claim 6:} $e_1, b_2, e_3, \ldots, e_d$ is a lattice basis.

\begin{proof}[Proof of Claim 6]

These elements form a basis of $\NR$ because of $e_2 = -e_1 -2 b_2$.
Let $e^*_1, b^*_2, e^*_3, \ldots, e^*_d$ denote the dual basis of $\MR$. By Figure \ref{picp2} and Equation (\ref{eq2}) we see
\[u_F = e^*_1 - b^*_2 + e^*_3 + \cdots + e^*_d,\]
\[u_{\nf{F}{e_2}} = u_F + b^*_2,\]
\[u_{\nf{F}{e_i}} = u_F - e^*_i \text{ for } i = 3, \ldots, d.\]
Since by reflexivity of $P$ the outer normals are lattice points in $M$, also $e^*_1, b^*_2, e^*_3,$\linebreak$\ldots, e^*_d$ are lattice points 
in $M$, thus $e_1, b_2, e_3, \ldots, e_d$ is a lattice basis of $N$.
\end{proof}

This finishes the proof of Proposition \ref{part3}, and hence of Theorem \ref{main}.
\end{proof}

\bibliographystyle{amsalpha}

\end{document}